\documentclass[draft]{amsart}
\usepackage{amssymb,latexsym,amsfonts,verbatim,amscd,ifthen}
\usepackage{color}

\numberwithin{equation}{section}

\theoremstyle{plain}

\newtheorem{theorem}[equation]{Theorem}   
\newtheorem{lemma}[equation]{Lemma}

\newtheorem{corollary}[equation]{Corollary}

\theoremstyle{definition}

\newtheorem{example}[equation]{Example}

\newtheorem*{maintheorem}{Main~Theorem}
\newtheorem*{mainproof}{Proof~of~Main~Theorem}

\DeclareMathOperator{\rank}{rank}

\begin{document}

\renewcommand{\:}{\! :\ }
\newcommand{\p}{\mathfrak p}
\newcommand{\m}{\mathfrak m}
\newcommand{\e}{\epsilon}
\newcommand{\g}{\gamma}
\newcommand{\lra}{\longrightarrow}
\newcommand{\ra}{\rightarrow}
\newcommand{\altref}[1]{{\upshape(\ref{#1})}}
\newcommand{\bfa}{\boldsymbol{a}}
\newcommand{\bfb}{\boldsymbol{b}}
\newcommand{\bfc}{\boldsymbol{c}}
\newcommand{\bfM}{\mathbf M}
\newcommand{\bfN}{\mathbf N}
\newcommand{\bfI}{\mathbf I}
\newcommand{\bfC}{\mathbf C}
\newcommand{\bfB}{\mathbf B}
\newcommand{\bsfC}{\bold{\mathsf C}}
\newcommand{\bsfT}{\bold{\mathsf T}}
\newcommand{\mc}{\mathcal}
\newcommand{\smsm}{\smallsetminus}
\newcommand{\ol}{\overline}
\newcommand{\twedge}
           {\smash{\overset{\mbox{}_{\circ}}
                           {\wedge}}\thinspace}

\newlength{\wdtha}
\newlength{\wdthb}
\newlength{\wdthc}
\newlength{\wdthd}

\title{UPPER BOUNDS FOR BETTI NUMBERS OF MULTIGRADED MODULES}
\author[Amanda Beecher]{Amanda Beecher}
\address{Department of Mathematics\\
         University at Albany, SUNY\\
         Albany, NY 12222}
\email{am2875@albany.edu}

\begin{abstract}
This paper gives a sharp upper bound for the Betti numbers of a finitely
generated multigraded $R$-module, where $R=\Bbbk [x_{1},\ldots,x_{m}]$ is the
polynomial
ring over a field $\Bbbk$ in $m$ variables.  The bound is given in terms of 
the rank and the first two Betti numbers of the module.  An example is given which 
achieves these bounds simultaneously in each homological 
degree.  Using Alexander duality, a bound is established for the total multigraded Bass numbers of 
a finite multigraded module in terms of the first two total multigraded Bass numbers.
\end{abstract}

\maketitle

\section*{Introduction}
Let $R=\Bbbk [x_{1},\ldots,x_{m}]$ be the polynomial ring over a field $\Bbbk$ in $m$ variables with the 
standard $\mathbb{Z}^m$ grading and $L$ a finite $\mathbb{Z}^m$
(multigraded) $R$-module.
Much work has been done on establishing lower bounds for
Betti numbers of $L$, initially motivated by the Buchsbaum-Eisenbud-Horrocks conjecture
for finite modules over regular local rings.  This conjecture was shown to hold for $R/I$ when $I$ is a monomial ideal
 by Evans and Griffith \cite{EG} and generally for all
multigraded modules by Charalambous \cite{Ch} and Santoni \cite{S}.

On the other hand, little is known about upper bounds for the Betti numbers.   
The main result of this paper gives sharp upper bounds for the total 
Betti number of $L$
 in each homological degree in terms of the rank and the first two 
Betti numbers of the module $L$.  
\begin{maintheorem}  For $i\geq 2$, we have $$\beta_{i}(L)\leq \binom{\beta_1(L)}{\beta _0(L)-\rank L+i-1}\binom{\beta_0(L)-\rank L+i-3}{i-2} .$$ 
\end{maintheorem} These bounds are precisely the ranks of the free modules in the multigraded Buchsbaum-Rim complex from \cite{CT} (called there the Buchsbaum-Rim-Taylor complex).  However, there seems to be no appropriate map that achieves comparison between the minimal free resolution of $L$ and the Buchsbaum-Rim complex.  Instead, we obtain our bound using the combinatorial structure of the (not necessarily minimal) free resolution defined by Tchernev \cite{T}.  We show that our bounds are sharp by giving  a class of examples that 
attains them simultaneously in each homological 
degree.  To achieve this, we make the generators of $L$ and their multidegrees sufficiently generic, so that the Buchsbaum-Rim-Taylor complex is the minimal free resolution.  In the last section, as a corollary to our main theorem, we give bounds for the total multigraded Bass numbers
 of multigraded modules in terms of the first two total multigraded Bass numbers, by using the Alexander duality 
functors defined by Miller \cite{M}.

\section{Counting T-flats}

A \emph{matroid} $\bfM$ is a finite set $S$ coupled with the nonempty
collection of \emph{independent}
subsets $\bfI$ of $S$ satisfying the following two properties:
\begin{enumerate}
\item
If $ Y\in \bfI \mbox{ and } X\subseteq Y$ then $X\in \bfI.$
\item
If $X,Y \in \bfI$, and $|Y|=|X|+1$ then there is $y\in Y\smsm X$
such that $X\cup \{y\} \in \bfI$.
\end{enumerate}
A maximal independent set is a \emph{basis}.  The collection of all bases
will be denoted by $\mathcal{B}(\bfM)$.
We say that a subset $A\subseteq S$ is \emph{dependent} if $A\not\in \bfI.$
A minimal dependent set is called a \emph{circuit}.  The collection of all
circuits of $\bfM$ will be denoted by $\mathcal{T}_{0}(\bfM)$.  Notice that
a matroid can be defined by specifying the set $\mathcal{B}(\bfM)$ or the
set $\mathcal{T}_{0}(\bfM)$.  The \emph{rank} in $\bfM$ of $A\subseteq S$ is
the number $r_{A}^{\bfM}= \max \{ |X| \mid X\in \bfI \mbox{ and } X\subseteq A\}$.  Notice
that $r_{X}^{\bfM}=|X|$ precisely when $X\in \bfI$ and $r_{A}^{\bfM}=|A|-1$
when $A\in \mathcal{T}_{0}(\bfM)$.  A subset $F\subseteq S$ is called a \emph{flat}
of the matroid $\bfM$ if it has the property that for every $x\in S\backslash F$ 
we have $r^{\bfM}_{F\cup \{x\}} =r^{\bfM}_{F} +1.$  The \emph{dual matroid} $\bf{M}^*$ is defined on 
the same finite set $S$ as $\bfM$ with $\mathcal{B}(\bfM ^*)=\{S\backslash B \mid B\in\mathcal{B}(\bfM)\}.$  
In the dual matroid, the rank
 of a set $A\subseteq S$ is the number $r^{\bfM^*}_A=|A|-r^{\bfM}_{S}+r^{\bfM}_{S\backslash A}.$

The subset $A$ is a \emph{T-flat} of $\bfM$ precisely when the complement $S\backslash A$ is 
a proper flat of the dual matroid $\bf{M}^*$ (\cite{T}, Definition 2.1.1).  The
\emph{level} of a subset $A\subseteq S$ is defined to be the number $$l_{A}=|A|-r_{A}^{\bfM}-1.$$  Notice that the level of a circuit is $0$, hence the notation
$\mathcal{T}_{0}(\bfM)$ for the collection of circuits.  We will similarly
denote the collection of T-flats of level $k$ by $\mathcal{T}_{k}(\bfM)$.  

\begin{lemma}\label{T:counting tflats}

Let $\bfM$ be a matroid on a finite set $S$.
If $ |S| =n$ and $r^{\bfM}_S = r$ then  $$|\mathcal{T}_{k}(\bfM)| \leq
\binom{n}{r+k+1}.$$
\end{lemma}
\begin{proof}
First, rewrite the rank of the dual matroid as follows
$$r^{\bfM^*}_A=|A|-r^{\bfM}_{S}+r^{\bfM}_{S\backslash A}=|S|-r^{\bfM}_S-l_{S\backslash A}^{\bfM}-1.$$
In the above notation, we see for a T-flat $A$ of level $k$, the rank of its complementary flat
 in the dual matroid will be 
$$r^{\bfM^*}_{S\backslash A}=n-r-k-1.$$
So by definition $$|\mathcal{T}_{k}(\bfM)|=|\mathcal{F}_{n-r-k-1}(\bfM^*)|$$ where $\mathcal{F}_{\rho}(\bfM)$ denotes the collection of all flats with rank $\rho$ in the matroid $\bfM$.
Since a subset has the property that $r_{A}^{\bfM}\leq |A|$, we see that the number of flats with rank $n-r-k-1$ 
must be less than the total number of subsets of $S$ with cardinality $n-r-k-1$.
Thus $$|\mathcal{T}_{k}(\bfM)|=|\mathcal{F}_{n-r-k-1}(\bfM^*)|\leq \binom{n}{n-r-k-1}=\binom{n}{r+k+1}$$ which verifies our bound for the number of T-flats of level $k$.
\end{proof}

\section{Betti Numbers of Multigraded Modules}

Let $$E\stackrel{\Phi}{\lra} G\lra L\lra 0$$ be a minimal finite free multigraded presentation of $L$.  We consider the field $\Bbbk$ as an $R$-module under the quotient map where
we send each variable to 1.  In this way, tensoring the free presentation with $\Bbbk$ gives the map
$$\begin{CD}
E \otimes_{R} \Bbbk 
@> \Phi \otimes id>>
 G\otimes_{R} \Bbbk =W
\end{CD} $$

Let $S=S\otimes 1$ be a basis of $E\otimes _{R} \Bbbk$.  Then the set map 
$$ 
\begin{CD}
\phi: S     
@> \Phi \otimes id >> 
W 
\end{CD} 
$$
defines a matroid $\bfM$ on $S$ where the independent sets are precisely those subsets $A$ of $S$ whose image under 
$\phi$ spans a vector space $V_A$ of dimension $|A|$.  
From this matroid Tchernev (\cite{T}, Section 2.2) constructs for each T-flat $I$ the vector 
space (called the T-space of $I$) $T_I(\phi)$.  Then the resolution of $L$ from (\cite{T}, Definition 4.3) 
has the form (with $\lambda=l_{S}+2=|S|-r_S^{\bfM}+1$), \[
T_\bullet(\Phi, S) =
0\rightarrow
T_\lambda(\Phi,S) \overset{\Phi_\lambda}\lra
T_{\lambda -1}(\Phi, S) \rightarrow \dots \rightarrow
T_1(\Phi, S) \overset{\Phi_1}\lra
T_0(\Phi, S) \rightarrow 0,
\] and the free $R$-modules are defined as
$$T_{i}(\Phi , S)=\bigoplus_{I\in \mathcal{T}_{i-2}(M)} R\otimes_{\Bbbk} T_I(\phi),\quad
2\leq i\leq l_S +2$$ and $$T_{0}(\Phi, S)=G \quad \mbox{ and }\quad
T_{1}(\Phi, S)=E.$$

\begin{mainproof}
Let $E\stackrel{\Phi}{\lra} G\lra L\lra 0$ be a minimal finite free multigraded presentation of $L$ so that 
$\beta_{0}(L) = \rank(G) \mbox{ and } \beta_{1}(L)= \rank(E).$
When $i\geq 2$ then $$\beta_{i}(L) \leq
\rank \bigl( T_{i}(\Phi,S)\bigr).$$
By definition
$$\rank \bigl(T_{i}(\Phi,S)\bigr) =\rank\left( \bigoplus_{I\in \mathcal{T}_{i-2}(\bfM)}
R\otimes T_I(\phi) \right) \leq |\mathcal{T}_{i-2}(\bfM)| \max_{I\in
\mathcal{T}_{i-2}(\bfM)}
(\dim_{\Bbbk} T_{I}(\phi)).$$ 
 It follows from the definition (\cite{T}, Definition 2.2.3) that for each $I \in \mathcal{T}_k(\bfM)$ one has $$\dim T_{I}(\phi) \leq \dim S_{k}(V)=\binom{r+k-1}{k},$$ where $S_k(V)$ denotes the $k^{th}$ symmetric power of $V=Im \phi$ and $r=\dim_{\Bbbk} (V)$.  Therefore, by Lemma \ref{T:counting tflats}, we get  
$$\beta_{i}(L) \leq \binom{|S|}{r^{\bfM}_{S}+i-1} \binom{r^{\bfM}_{S}+i-3}{i-2}
.$$  Since $$|S|=\rank (E)=\beta _1(L)$$ and 
$$r^{\bfM}_S =\rank (\phi)=\rank (\Phi)=\rank G-\rank L=\beta _0(L)-\rank L,$$ the bounds on Betti numbers of
multigraded modules have been established.
\end{mainproof}

Next, we give an example of a finite multigraded module that achieves these bounds simultaneously in each homological degree.

\begin{example}
Let $L$ be a finite 
multigraded module with a minimal presentation matrix of uniform rank.  In the above set up, this is equivalent to saying that if the rank of the matrix $\phi$ is $r$ then the image of every $r$-element subset of $S$ has dimension $r$.  Then the free $R$-modules described in \cite{T} (Definition 4.2) are precisely those of the Buchsbaum-Rim-Taylor complex 
defined in \cite{CT} (Definition 4.3).  In addition, choose the multidegrees of the 
generators of $E$ to be generic.  In that case no two T-flats have the same multidegree.  
Therefore, the Buchsbaum-Rim-Taylor complex is the minimal free multigraded resolution of $L$.  
Thus for $i\geq 2$, $$\beta_{i} (L)=\binom{|S|}{r+i-1}\binom{r+i-3}{i-2} ;$$ our bounds are simultaneously achieved.
\end{example}

\section{Multigraded Bass numbers of multigraded modules}

We write $\mu _{i}(\p ,L)$ for the $i^{th}$ Bass number of an $R$-module $L$
at a prime $\p$.  Since we are over a polynomial ring and $L$ is 
multigraded, there are only finitely many mulitgraded primes $\p$ generated by $\{x_{j_{1}},\ldots,x_{j_{l}}\} \text{ for } \{j_1,\ldots,j_l\} \subset \{1,\ldots, m\}.$  Thus it makes sense to call the well-defined integer $$ ^*\mu _i(L)=\sum_{\p: \text{multigraded prime ideal}} \mu_i(\p,L)$$ the total multigraded $i^{th}$ Bass number of $L$.

\begin{theorem}\label{T:BassBounds}
Let
$R=\Bbbk[x_1,\dots, x_m]$ be a polynomial ring over
a field $\Bbbk$ with the standard grading,
$L$ a finite multigraded $R$-module, then the bounds for $i\geq 2$ of the total multigraded Bass numbers of $L$ are
$$ ^*\mu_{i}(L)\leq  \binom{^*\mu _1(L)}{^*\mu_0(L)+i-1}\binom{^*\mu _0(L)+i-3}{i-2}
.$$
\end{theorem}

\begin{proof}
For any $\bfc$, $^*\mu_i (L(-\bfc))=\ ^*\mu _i(L)$.  Since $L$ is finitely generated and has a finite multigraded injective resolution, there is a $\bfc$ so that $L(-\bfc)$ only has non-zero degrees greater than or equal to $\mathbf{1}$ and so that the nonzero Bass numbers at each multigraded prime occur in only positive degree.  Thus it suffices to establish our bound when $L$ satisfies these two conditions.  

Let $\bfa$ be the componentwise maximum of the degrees of the $0^{th}$ and $1^{st}$ Betti numbers of $L$.  Thus $L$ is a positively $\bfa$-determined module as
 defined in (\cite{M}, Definition 2.1) and the minimal free resolution will be positively $\bfa$-determined (\cite{M}, Proposition 2.5). 
 
 We
write $L^{\bfa}$ for the Alexander dual of an $R$-module $L$ with respect to degree $\bfa$ as defined in \cite{M}.
We will write $B_{\bfa} L$ of an $R$-module for the quotient of $L$ by the submodule $\bigoplus _{\bfb\nleq \bfa}L_{\bfb}$ as in \cite{M}.  In the case when $L$ is generated in degrees greater than $\mathbf{1}$, we have that $B_{\bfa} L^{\bfa}=L^{\bfa}$, since $L^{\bfa}$ is bounded inside the interval $[\mathbf{0}, \bfa-\mathbf{1}]$.  Thus we will use Miller's results for $B_{\bfa} L^{\bfa}$ and $L^{\bfa}$ interchangeably throughout the remainder of this paper.  

Let $$0\lra  L\lra I \overset{\Lambda}\lra J $$ be a minimal multigraded copresentation of $L$.  Miller (\cite{M}, 
Theorem 4.5) shows that the matrix $\Lambda$ is also a minimal presentation matrix for a free resolution 
of $L^{\bfa}$ after some appropriate degree shifts.  In establishing bounds for the Betti numbers 
in the previous section, we obtain the coefficient matrix of the vector space complex, which is independent 
of the degrees of the original free modules.  
Thus the coefficient matrix 
of $\Lambda$ defines a representable matroid for the multigraded module $L^{\bfa}$.  
One easily sees that the bounds established in the Main 
Theorem depend only on the presentation matrix.  Since $\rank L^{\bfa}=0$ we have for $i\geq 2$
$$\beta_{i}(L^{\bfa}) \leq  \binom{\beta_1 (L^{\bfa})}{\beta _0(L^{\bfa})+i-1}\binom{\beta _0(L^{\bfa})+i-3}{i-2}$$

Further, Miller (\cite{M}, Theorem 5.3) establishes the relation 
$$\beta_{i,\bfb}(L^{\bfa})=\mu_{i,(\bfa \backslash \bfb)-supp(\bfb)} (\m^{supp(\bfb)},L)$$
where $0\leq \bfb \leq \bfa \cdot supp(\bfb)$.
By summing over all possible degrees and since the nonzero Bass numbers at each multigraded prime occur in positive degrees, we see that $^*\mu _0(L)=\beta_0 (L^{\bfa})=$ number of rows of $\Lambda$, $^*\mu _1(L)=\beta_1 (L^{\bfa})=$ number of columns of $\Lambda$, and 
$$^*\mu_{i} (L)=\beta_{i}(L^{\bfa}).$$
Therefore,
the bounds
for the total multigraded Bass numbers have been established.
\end{proof}

The following corollary generalizes this bound to all Bass numbers.

\begin{corollary}
Let
$R=\Bbbk[x_1,\dots, x_m]$ be a polynomial ring over
a field $\Bbbk$ with the standard $\mathbb{Z}^m$ grading,
$L$ a finite multigraded $R$-module and $\p$ a prime ideal of $R$.  Let $d=\dim R_{\p}/\p ^* R_{\p}$, where $\p ^*$ denotes the largest multigraded prime ideal of $R$ contained in $\p$.  Then, for $i\geq 2+d$, the bounds for Bass numbers of $L$ are
$$ \mu_{i}( \p , L)\leq  \binom{^*\mu _1(L)}{^*\mu_0(L)+i-d-1}\binom{^*\mu _0(L)+i-d-3}{i-d-2}
.$$
\end{corollary}

\begin{proof}
Goto and Watanabe showed in \cite{GW} that $$\mu _i(\p , L)=\mu _{i-d}(\p ^*, L).$$  
Clearly, $\mu_{i-d}(\p ^*, L) \leq {} ^*\mu_{i-d}(L)$.  Thus Theorem \ref{T:BassBounds} establishes our bounds.
\end{proof}

\end{document}